\magnification=\magstep1
\documentstyle{amsppt}
\define \res{\operatornamewithlimits{res}}
\hsize=6.5 true in
\vsize=9 true in
\topmatter
\title   High Moments of the Riemann Zeta--Function
\endtitle
\author  J. B. Conrey and S. M. Gonek
\endauthor
\address
J. B. Conrey\endgraf
American Institute of Mathematics\endgraf
360 Portage Ave.\endgraf
Palo Alto, CA 94306\endgraf
{\it E-mail address:}  {\bf conrey\@aimath.org\endgraf}
\null
Department of Mathematics\endgraf
Oklahoma State University\endgraf
Stillwater, OK 74078-0613\endgraf
\null
S. M. Gonek\endgraf
Department of Mathematics \endgraf
University of Rochester\endgraf
Rochester, NY 14627\endgraf
{\it E-mail address:}  {\bf gonek\@math.rochester.edu}
\endaddress
\thanks Research of both authors was supported in part by the
American Institute of Mathematics and by grants from the NSF.
\endthanks
\endtopmatter

\def\D{\Bbb D}
\def\E{\Bbb E}

\NoBlackBoxes

\centerline{\bf Introduction}

\medskip

One of the most important goals of number theorists this
century has been to determine the moments of the Riemann
zeta--function on the critical line. These are important
because they can be
used to estimate the maximal order of the zeta--function
on the critical line, and because of their applicability
to the study of the distribution of prime numbers, often
through zero--density estimates, and to divisor problems.

The two most significant early results were obtained by
Hardy and Littlewood [HL] in 1918 and Ingham [I] in 1926.
Hardy and Littlewood proved that
$$
\int_0^T|\zeta(1/2+it)|^2~dt \sim T \log T \tag 1
$$
as $T\to \infty$, and Ingham showed that
$$
\int_0^T|\zeta(1/2+it)|^4~dt \sim \frac {1}{2\pi^2} T \log ^4 T \ . \tag 2
$$

No analogous formula has yet been proved for any higher moment,
and it seems unlikely that any will be in the near future.
In fact, the problem is so intractable that, until a few years ago,
no one was even able to produce a plausible guess for the
asymptotic main term. Recently, however, Conrey and Ghosh [CG2]
found a special argument in the case of the sixth power moment that
led them to conjecture that
$$
\int_0^T |\zeta(1/2+it)|^6~dt \sim \frac {42}{9!}\prod_p\left(\left(1-\frac
1p\right)^4
\left(1+\frac 4 p+\frac 1{p^2}\right)\right)T\log^9 T \ . \tag 3
$$

The object of this paper is to describe a new heuristic approach
that leads to a conjecture for the asymptotic main term of
$$
I_k(T)=\int_0^T|\zeta(1/2+it)|^{2k}~dt \ ,
$$
when $k$ equals three and four.
The resulting formula for the sixth power moment is identical
to the conjecture of Conrey and Ghosh above
and lends additional strong support to it in a sense to be
described below. For the eighth power moment, we obtain the
following new conjecture.

\proclaim {Conjecture 1} As $T \to \infty$,
$$
\int_0^T |\zeta(1/2+it)|^8~dt \sim
\frac{24024}{16!}\prod_p\left(\left(1-\frac 1p\right)^{9}
\left(1+\frac9p+\frac9{p^2}+\frac1{p^3}\right)\right)T\log^{16}T.
$$
\endproclaim

We will also discuss the size of higher moments of the
zeta--function and its maximal order in the critical strip.

While writing this paper, the authors learned
that J. Keating and N. Snaith [KS] have made a high moments
conjecture based on a completely different approach.
Instead of the attack through approximate functional
equations, mean value theorems, and additive divisor sums
employed here, they prove a general result on
moments of random matrices whose eigenvalues have a
GUE (Gaussian Unitary Ensemble) distribution. If the zeta--
function is modeled by the determinant of such
a matrix, and there are reasons to believe it is,
then the moments they calculate apply to the zeta--function
as well. It is remarkable that our conjecture and theirs,
which we state later, agree for the sixth and eighth moments,
and it suggests that both are likely to be right.

We begin by outlining the main ideas behind our approach, 
starting with a brief discussion of approximate functional equations.

For $s=\sigma+it$ and $\sigma>1$, $\zeta^k(s)$ has the
Dirichlet series expansion
$$
\zeta^k(s)=\sum_{n=1}^\infty\frac{d_k(n)}{n^s} \ ,
$$
where $d_k(n)$ is the $k$th divisor function, which
is multiplicative and defined at prime powers by
$d_k(p^j)=\binom{k+j-1}{j}$.
The series does not converge when
$\sigma\le1$, but we can nevertheless approximate
$\zeta^k(s)$ in this region by a sum of two Dirichlet
polynomials. This is called an approximate functional
equation, and its prototype is
$$
\zeta(s)^k=\D_{k,N}(s)+\chi(s)^k \D_{k,M}(1-s)+\E_k(s)\ , \tag 4
$$
where
$$
\D_{k,N}(s)=\sum_{n=1}^N \frac{d_k(n)}{n^s}\ \ ,
$$
$\E_k(s)$ is an error term, $MN=\left(\frac{t}{2\pi}\right)^k$,
and
$$
\chi(s)=(\pi)^{s-1/2}\frac{\Gamma(\frac{1-s}{2})}{\Gamma(\frac{s}{2})}
$$
is the factor from the functional equation for the zeta--function,
namely
$$
\zeta(s)=\chi(s)\zeta(1-s)\ .
$$
Note that from the last equation it follows that $\chi(s)$ satisfies
$$
\chi(s)\chi(1-s)=1\ .
$$

Taking $s=1/2+it$ in (4), integrating the square of
the modulus of both sides, and assuming that $\E_k(1/2+it)$ is
sufficiently small, we obtain
$$
\split
\int_{T}^{2T}|\zeta(1/2+it)|^{2k}~dt
&\sim
\int_{T}^{2T}|\D_{k,N}(1/2+it)|^2~dt+\int_{T}^{2T}|\D_{k,M}(1/2+it)|^2~dt\\
&+2\Re\int_{T}^{2T}\chi(1/2-it)^k\D_{k,N}(1/2+it)\D_{k,M}(1/2+it)~dt\ .
\endsplit  \tag 5
$$
Now
$$
\chi(1/2-it)=\exp\left(it\log \frac{t}{2\pi e}\right)(1+O(1/t))
$$
as $t\to \infty$, so we find that
$$
\chi(1/2-it)^k(mn)^{-it}=\exp\left(it\log \frac{(t/2\pi e)^k}
{mn}\right)(1+O(1/t))\ .
$$

This has a stationary phase at $t=2\pi(mn)^{1/k}$, which
is generally outside the interval of integration. This
suggests that when $ MN \leq T^{k - \epsilon} $,
the third integral on the right--hand side
of (5) is smaller than the larger of the first two.
(That it is no larger can be seen from the Cauchy--Schwarz inequality.)
An extrapolation from this to
 $ MN = T^{k} $ seems reasonable, and
one can probably also replace $t$ in the condition
$MN=\left(\frac{t}{2\pi}\right)^k$  by~$T$ when $t$ is
large. One can also show that if the Lindel\"{o}f Hypothesis
is true, and if $M=0$ and $N \gg T^{k}$, then
$\zeta(1/2+it)^{k}$ is well approximable in mean square by
$\sum_{n=1}^{\infty}  {d_k(n)}e^{-n/N}{n^{-s}}$.
This suggests that
$\zeta(1/2+it)^{k}$ should also be well approximable
by $ \D_{k,N}(1/2+it) $  for the same $M$ and $N$\.
Thus, we expect the following to hold.

\proclaim {Conjecture 2}
For every positive integer~$k$, we have
$$
\int_{T}^{2T}|\zeta(1/2+it)|^{2k}~dt
\sim \int_{T}^{2T}|\D_{k,N}(1/2+it)|^2~dt+\int_{T}^{2T}
|\D_{k,M}(1/2+it)|^2~dt, \tag 6
$$
where
$$
MN=\left(\frac{T}{2\pi}\right)^k
\
\hbox{with}
\
M,N\ge 1/2,\ \ \ \ \hbox{or}\ \ \ \ N\gg\left(\frac{T}{2\pi}\right)^k
\hbox{if}\ M=0. \tag 7
$$
\endproclaim

By classical methods one can prove that Conjecture 2
holds when $k=1$, and also when $k=2$ provided that
$\max(M,N)\ll T$. When $k\ge3$, however, the known bounds
for $\E_k(s)$ in (4) are too large to give (6), and it is
also difficult to show that the third term in (5) really
is smaller than the other two. Nevertheless, it may be possible
to overcome these problems (when $k= 3$ or $4$) by appealing
to a more complicated
form of the approximate functional equation first developed
by A. Good [Go] for $\zeta(s)$ and for the $L$--functions
attached to cusp forms (which are analogous to $\zeta(s)^2$).
We hope to return to this question in a future article.

Our problem now reduces to determining an asymptotic
estimate for the mean square of the Dirichlet polynomial
$\D_{k,N}(1/2+it)$. The standard tool for this is the
classical mean value theorem for Dirichlet polynomials which,
in the refined version due to Montgomery and Vaughan [MV],
asserts that
$$
\int_{T}^{2T}\left|\sum_{n=1}^N a(n)
n^{it}\right|^2~dt=\sum_{n=1}^N|a(n)|^2(T+O(n)) \,. \tag 8
$$
Using this, we see that
$$
\int_{T}^{2T}\left|\D_{k,N}(1/2+it)\right|^2~dt=
\sum_{n\le N}\frac{d_k(n)^2}{n}(T+O(n))\ .
$$
Now it is well known that
$$
\sum_{n\le N} d_k(n)^2 \sim
\frac{a_k}{\Gamma(k^2)}N \log^{k^2-1}N
$$
and
$$
\sum_{n\le N} \frac {d_k(n)^2}{n} \sim
\frac{a_k}{\Gamma(k^2+1)} \log^{k^2}N\ , \tag 9
$$
where
$$
a_k=\prod_p\left(\left(1-\frac {1}{p}\right)^{k^2}\sum_{r=0}^{\infty}
\frac{d^{2}_{k}(p^{r})}{p^r}\right)\ . \tag 10
$$
Thus we deduce that
$$
\int_{T}^{2T}\left|\D_{k,N}(1/2+it)\right|^2~dt \sim
\frac{a_k}{\Gamma(k^2+1)}T \log^{k^2}N\ ,
$$
for $N\ll T$, say. Inserting this into (6), and assuming that
$M$ and $N$ satisfy (7) and are each $\ll T$,
we arrive at the previously known estimates for $I_1(T)$
and $I_2(T)$ in (1) and (2).

The problem we encounter for higher moments is that
at least one of $M$ and $N$ must be significantly
larger than $T$ because of (7), but then the $O$-terms
dominate the right--hand side of (8) and we lose the
asymptotic formula. To get around this we appeal to recent
work of Goldston and Gonek [GG] on mean values of ``long''
Dirichlet polynomials, which allows one to evaluate
$$
\int_{T}^{2T} \left|\sum_{n=1}^N \frac{a(n)}{n^{1/2+it}}\right|^2~dt \tag 11
$$
provided one has a good handle on the coefficient sum
$$
A(x)=\sum_{n\le x} a(n)
$$
and on the coefficient--correlation sums
$$
 A(x,h)=\sum_{n\le x} a(n) a(n+h)\ .
$$
Specifically, one needs formulae of the type
$$
A(x)=m(x)+E(x)
$$
and
$$
A(x,h)=m(x,h)+E(x,h)
$$
in which $m(x)$ and $m(x,h)$ are differentiable with respect to $x$,
$E(x)\ll x^{\theta}$, and $E(x,h)\ll x^{\phi}$, uniformly for
$h\ll x^{\eta}$ with $0< \theta, \phi, \eta <1$.
One can then obtain an asymptotic formula for (11)
for $N \ll T^{\min(\frac{1}{\theta},
\frac{1}{\phi}, \frac{1}{1-\eta})-\epsilon}$
and any $\epsilon>0$. If, in addition, one knows something
about averages of the error terms $E(x,h)$ with respect to $h$,
$N$ can be taken even larger. For example (confer [GG]), taking
$a(n)=d_k(n)$ and $k=1$, we see that we can choose $\theta$ and
$\phi$ to be $\epsilon$ and $\eta$ to be $1-\epsilon$ for
any small positive $\epsilon$, which means that $N$ can be an
arbitrarily large power of $T$.
When $k=2$, a result of Heath-Brown [H-B1] allows us
to take $\phi=\eta=5/6$, and
this permits us to take $N$ up to $T^{6/5-\epsilon}$.
(A better result can probably be obtained by using more
recent tools, such as those of [DFI], combined with an averaging over h.)
However, what we actually expect in this case is that $\phi=1/2+\epsilon$
and $\eta\ge 1/2-\epsilon$, and if this is the case,
then our method leads to Ingham's formula (2) for all $N\ll T^{2-\epsilon}$.

For $k\ge 3$, unfortunately, the estimates we require
for the additive divisor sums
$$
D_k(x,h)=\sum_{n\le x}d_k(n)d_k(n+h)
$$
have never been proved. In fact, even an asymptotic formula for
$$
\sum_{n\le x}d_3(n)d_3(n+1)\ ,
$$
is not known. Still, a precise formula for the main term
of $D_k(x,h)$ can be conjectured
in several ways, the easiest being the so-called
$\delta$--method of Duke, Friedlander, and Iwaniec [DFI].
In the next section we describe how their method,
together with a guess as to how the error term behaves,
suggests

\proclaim{Conjecture 3} Let
$D_k(x,h)=\sum_{n\le x}d_k(n)d_k(n+h)$.
Then we have
$$
D_k(x,h)=m_k(x,h)+O(x^{1/2+\epsilon})\tag 12
$$
uniformly for $1\le h \le x^{1/2}$, where $m_k(x,h)$ is
a smooth function of $x$. The derivative of $m_k(x,h)$ is given by
$$
m_{k}'(x,h)=\sum_{d\mid h} \frac{f_k(x,d)}{d}\ , \tag 13
$$
where
$$
f_k(x,d)=\sum_{q=1}^\infty \frac{\mu(q)}{q^2}P_k(x,qd)^2\ , \tag 14
$$

$$
P_k(x,q)=\frac{1}{2\pi i}\int_{|s|=1/8}
\zeta^k(s+1)G_k(s+1,q)\left(\frac{x}{q}\right)^s~ds\ , \tag 15
$$

$$
G_k(s,q)=\sum_{d\mid q}\frac{\mu(d)}{\phi(d)}d^s\sum_{e\mid d}
\frac{\mu(e)}{e^s}g_k(s,qe/d)\ , \tag 16
$$
and, if $q=\prod_pp^{\alpha}$,
$$
g_k(s,q)=\prod_{p\mid
q}\left(\left(1-p^{-s}\right)^k\sum_{j=0}^\infty
\frac{d_k\left(p^{j+\alpha}\right)}{p^{js}}\right)\ . \tag 17
$$
Furthermore, for $d\le x$ we have
$$
f_k(x,d)\ll d^{2}_{k-1}(d) \log^{2k-2}x\ . \tag 18
$$
\endproclaim

One can actually prove the last assertion and we shall do so below.
A similar conjecture for $m_{k}(x,h)$ appears in Ivic [Iv],
but its form is less appropriate for our purposes.

Using Conjecture 3 together with Theorem 1 of Goldston
and Gonek ~[GG], we are led to

\proclaim{Conjecture 4} Let $N=T^{1+\eta}$ with $0\le \eta \le 1$.
Then
$$
\int_T^{2T}\left|\D_{k,N}(1/2+it)\right|^2~dt\sim w_{k}(\eta)
\frac{a_k}{\Gamma(k^2+1)}TL^{k^2}\ ,
$$
where $a_k$ is given by (10), and
$$
w_{k}(\eta)=(1+\eta)^{k^2}\left(1-\sum_{n=0}^{k^2-1}\binom{k^2}{n+1}
\gamma_k(n) \left(1-(1+\eta)^{-(n+1)}\right)\right)\ ,
$$
where
$$
\gamma_k(n)=(-1)^n\sum_{i=0}^k\sum_{j=0}^k\binom{k}{i}
\binom{k}{j}\binom{n-1}{i-1,j-1,n-i-j+1}
$$
for $n\ge 1$, and
$$
\gamma_k(0)=k.
$$

\endproclaim

\medskip

With more work, we could {\it prove} that Conjecture 3 implies
the mean value formula of Conjecture 4 for $0 \leq \eta \leq 1 -
\epsilon.$  However, we would then still have to extrapolate
to the full interval
$0 \leq \eta \leq 1 $. Originally we expected the error
terms in Conjecture 3 to exhibit considerable cancelation
when summed over h ranges of length up to $x^{1-\epsilon}$, just
as we believe they do when the coefficients $d_k(n)$ are
replaced by $\Lambda(n)$ or $\mu(n)$.
We were surprised to find, however, that this does not seem to be
the case. For example, when $k=2$, it appears that the error term in
(12) actually leads to a main term in (2) if we use (6) with $N\gg
T^{2+\epsilon}$. Thus, in the formula
$$
w_2(\eta)=1+4 \eta -6 \eta^2+4 \eta^3-\eta^4,
$$
if we take $\eta>1$, we see that $w_2(\eta) > 2$.
It follows that, (6), (7), and Conjecture 3  cannot all be
true if $\eta > 1$,
since they contradict Ingham's result (2).
We believe (6) and (7) are correct, so our conclusion is that
Conjecture 3 fails when $\eta >1$.  Since we have accounted
for all of the obvious main terms, we are forced to conclude
that somehow the error terms in (12) accumulate to deliver
a new ``mysterious'' main-term.

The restriction on the size of $N$ in Corollary 4 means
that $M$ and $N$ must satisfy
$\max(M,N)\ll T^2$ in (6) and,
in light of (7), this forces $k$ to be less than or equal to $4$.
When $k=3$, for example, we may use (6)with $N$ satisfying
$T\ll N \ll T^2$.
Then $M=T^3/N$ will satisfy the same bound and we will show
that for any choice of $N$ in this range, Conjecture 4 leads to
$$
\eqalign{
\int_{T}^{2T}|\zeta(1/2+it)|^{6}~dt &\sim
\int_{T}^{2T}|\D_{3,N}(1/2+it)|^2~dt+
\int_{T}^{2T}|\D_{3,T^{3}/N}(1/2+it)|^2~dt\cr
&\sim 42\frac{a_3}{9!}T\log^9 T\ .
}
$$
Adding this up for $T$ replaced by $T/2, T/4, \dots $,
we obtain
$$
I_{3}(T)\sim 42\frac{a_3}{9!}T\log^9 T\ \,,
$$
which is (3). The persistence of this estimate throughout
the range $T\ll N \ll T^2$ gives very strong independent
confirmation of the sixth power moment
conjecture of Conrey and Ghosh.
When $k=4$ we are forced to take $N=M=T^2$ and,
in this case, Conjecture 4 leads to
$$
\eqalign{
\int_{T}^{2T}|\zeta(1/2+it)|^{8}~dt
\sim 2\int_T^{2T}|\D_{4,N}(1/2+it)|^2~dt
\sim 24024\frac{a_4}{16!}T\log^{16}T\ ,
}
$$
so that
$$
I_{4}(T)\sim 24024\frac{a_4}{16!}T\log^{16}T\ .
$$

The form of these results suggests that
there exists a constant $g_k$
such that
$$
I_k(T)\sim g_k \frac{a_k}{\Gamma(1+k^2)}T\log ^{k^2} T \ . \tag 19
$$
With this notation, Hardy and Littlewood's result asserts that $g_1=1$,
Ingham's that $g_2=2$, the conjecture of
Conrey and Ghosh that $g_3=42$, and Conjecture 1
that $g_4=24024$. Since we do not know whether
$g_k$ exists in general, it is convenient to define
$$
g_k(T)=\left(\frac{a_k}{\Gamma(1+k^2)}T\log^{k^2}T\right)^{-1}I_k(T)\ ,
$$
so that
$$
g_k=\lim_{T \to \infty}g_k(T)\ ,
$$
provided the limit exists.

There have been numerous papers devoted to the estimation
of $g_k(T)$. For example, writing $g_k(T)\succeq C$ to mean that
$g_k(T)\ge (1+o(1))C$,
Conrey and Ghosh [CG3] showed unconditionally that $g_3(T)\succeq 10.13$
as $T\to \infty$. Soundararajan [S] increased the bound to $20.26$
and later (unpublished) to $24.59$. Subject to the truth of the
Lindel\"{o}f Hypothesis, Conrey and Ghosh [CG3] also obtained the
lower asymptotic bounds $g_4(T) \succeq 205$,
$g_5(T) \succeq 3242$, and $g_6(T)\succeq 28130$.
Their method uses the auxiliary means
$$
\int_0^T|\zeta(1/2+it)|^2\D_{k,N}(1/2+it)^2~dt \tag 20
$$
and
$$
\int_0^T|\zeta(1/2+it)|^2\zeta(1/2+it)^k \D_{k,N}(1/2-it)~dt \tag 21
$$
for $N=T^\theta$ with $0<\theta<1/2$.
Under the assumption that the Lindel\"{o}f Hypothesis is true and that
these formulae
also hold when $\theta$ tends to $1$, Conrey and Ghosh deduced the stronger
lower bounds
$g_3(T) \succeq 38.76$, $g_4(T)\succeq 21528$, $g_5(T)\succeq 48438800$, and
$$
g_k(T)\succeq (ek/2)^{2k-2}\ , \tag 22
$$
as $k\to \infty$.
Their conjecture that $g_3=42$ was based on similar ideas.

The function $g_k(T)$ has also been studied for non integral values of k
by Ramachandra, Ramachandra and Balasubramanian, Heath-Brown,
Conrey and Ghosh, Gonek, and Soundararajan, among others.
Summarizing just a few of the results:
Ramachandra [R] proved that $g_k(T)\approx 1 $ for $k=1/2$.
Heath-Brown [H-B2]
extended this to $k=1/n$ for $n$ any positive integer.
Conrey and Ghosh [CG1]
showed that the Riemann Hypothesis (RH) implies $g_k(T)\ge 1$ for all $k>0$
and
Soundararajan [S] improved this for all $k\ge 2$ by showing that $g_k(T)\ge 2$.
Gonek [G] proved that on RH $g_k(T)\ge 1$
for $-1/2 < k < 0$. Conrey and Ghosh [CG3] also gave conjectural
improvements in the lower bound for the interval $1<k<2$
assuming their conjecture that $\theta =1$ is permissible in (20) and
(21).


One rationale for studying $g_k(T)$ for non integral $k$
is that if $g_k$ exists and is meromorphic as a function
of $k$, then it can be identified from its values for
small real $k$. In fact, the other components in the
conjectural formula (19) for $I_k(T)$ are known to be entire
functions of order 2. This is clearly the case for
$\log^{k^2}T$ and $1/\Gamma(1+k^2)$, and was proved for $a_k$
by Conrey and Ghosh [CG3]. It would therefore be
interesting to know how $g_{k}$ behaves as $k\to \infty$ .
The conjectural result (22) suggests that $g_k$ grows
at least like a function of order 1.
We will see below that the function $w_{k}(\eta)$
in Conjecture 4 is $\succeq (1+\eta)^{k^2}$, so that
$$
\int_T^{2T}\left|\D_{k,T^2}(1/2+it)\right|^2~dt
\succeq {2^{k^2}} \frac{a_k}{\Gamma(1+k^2)}T\log^{k^2}T \ .
$$
Using this in (6) and dropping the second term,
which is positive and probably much larger than
the first, we deduce that
$g_k(T)\succeq 2^{k^2}$
as $k \to \infty$ through the integers.
Thus, Conjecture 4 implies that if $g_k$ exists
it grows at least as fast as a function of order 2.
Probably it grows no faster than this. To see why, take $M=N$
in (6) to obtain
$$
\int_T^{2T}|\zeta(1/2+it)|^{2k}~dt\sim 2\int_T^{2T}|\D_{k,N}(1/2+it)|^2~dt
$$
with $N=T^{k/2}$. According to (9), the contribution to this
of the ``diagonal'' terms is
$$
2T\sum_{n\le N}\frac{d_k(n)^2}{n}\sim
\frac{2 a_k}{\Gamma(k^2+1)}T \log^{k^2}N \ .
$$
But we expect this to be larger than the
entire mean value, as it is when $N\ll T$ by Montgomery
and Vaughan's mean value theorem, and when $T\ll N\ll T^2$
by Conjecture 4. This reasoning suggests that
$g_k\le 2(k/2)^{k^2}$. Thus, we believe that
$$
2^{k^2}\preceq g_k\preceq 2(k/2)^{k^2}\ .
$$
This is consistent with the conjecture of Keating
and Snaith referred to above, which is that
$$
g_k=\Gamma(1+k^2)\lim_{N\to
\infty}N^{-k^2}\prod_{j=1}^N\frac{\Gamma(j)\Gamma(j+2k)}{\Gamma(j+k)^2}.
$$
By Stirling's formula we easily see that this implies
$$
g_k=(k/4e^{1/2})^{k^2(1+o(1))}  \tag 23
$$
as $k\to \infty$.

With conjectural estimates for $g_k$ in hand we can approach
the question of the maximal order of the zeta--function.
Define
$$
m_T=\max_{0\le t \le T}|\zeta(1/2+it)|\ ,
$$
and for convenience write
$$
L=\log T.
$$
On RH it is known that
$$
m_T\ll \exp\left(C_u\frac{L}{\log L}\right) \tag 24
$$
for some positive constant $C_u$. On the other hand,
it follows from work of Montgomery [M] that if RH is true, then
$$
m_T \gg \exp \left(C_l \sqrt{\frac{L}{\log L}}\right) \tag 25
$$
with $C_l =1/20$. Subsequently Balasubamanian and
Ramachandra [BR] eliminated the need for RH in
the lower bound and Balasubramanian [B] increased the constant to
$C_l=0.5305... $. (The constant is quoted as ``3/4'' in his paper, but K.
Soundararajan has pointed out
that there is an error in the computation of $\max D(\ell)$ there; it seems
to be larger by a factor of $\sqrt{2}$ than
it should be.
The wide disparity between the upper and lower bounds here
appears in several other problems as well, and for the same reasons.
For example, on RH it is known that
$S(T) = \frac{1}{\pi} \arg \zeta(1/2+iT)$ satisfies
$$
S(T)\ll L/\log L\ , \tag 26
$$
and also (see [M]) that there exists a sequence of
values of $T\to \infty$ such that
$$
S(T) \gg \sqrt{\frac{L}{\log L}}\ .
$$
On the 1-line, the disparity appears as a factor of 2. Namely, RH implies that
$$
|\zeta(1+it)|\preceq 2 e^{\gamma} \log \log t \tag 27
$$
as $t \to \infty$ while, unconditionally, there exists a sequence of
$t\to \infty$ for which
$$
\zeta(1+it) \succeq e^{\gamma} \log \log t\ .
$$
The $q$-analogue of this asserts that if the Generalized
Riemann Hypothesis is true, then
$$
|L(1,\chi)|\preceq 2 e^{\gamma} \log \log q  \tag 28
$$
for every primitive character $\chi \pmod{q}$, whereas
unconditionally there is a sequence of $q \to \infty$
such that
$$
L(1,\chi_q)\succeq  e^{\gamma} \log \log q\ ,
$$
with $\chi_q$ a quadratic, primitive character $\pmod{q}$.
(See Shanks [Sh] for a discussion of his extensive
numerical work on this question.)

We can obtain lower bounds for $m_T$ directly from lower
bounds for $I_k(T)$ by observing that
$$
 m_T\ge
 \left(\frac{1}{T}\int_0^T|\zeta(1/2+it)|^{2k}~dt\right)^{1/2k}\ .
$$

To estimate the right--hand side of this we require an estimate for
$a_k$ in addition to our conjectural estimates for $g_k$.
To this end we shall prove the following

\proclaim{Proposition} Let
$a_k$ be defined as above. Then we have
$$
\log a_k =-k^2\log(2e^\gamma\log k)+o(k^2)
$$
as $k\to \infty$.
\endproclaim

With more work, and assuming RH, we could obtain the much
more precise result
$$\eqalign{
\log a_k =&-k^2 \left(\log \log k + \log(2e^{\gamma}) +
\log(1+ \log 2/ \log k)\right)\cr
&+8 k^2\left(\int_1^\infty\frac{\log(J_0(iw))}{w^3\log(4k^2/w^2)}~dw
+\int_0^1\frac{\log(J_0(iw))-w^2/4}{w^3\log(4k^2/w^2)}~dw \right)\cr
&+O(k^{1+\epsilon}),
}$$
where $J_0$ is the Bessel function of the first kind of order $0$. The
integrals can then be
expanded further to give an asymptotic expansion in decreasing
powers of $\log k$.

We now calculate a lower bound for $m_T$ assuming a lower bound for
$g_k$ of the form $(Ak+B)^{k^2}$.
By Stirling's formula we have
$$
m_T\succeq \left(\frac{g_k a_k L^{k^2}}{\Gamma(1+k^2)}\right)^{1/2k}
\succeq
\left(\frac{(Ak+B)Le^{1-\gamma}}{2k^2\log k}\right)^{k/2}\ . \tag 29
$$

It is not difficult to see that when $A=0$, the right-hand
side is maximized (as a function of k) essentially by taking
$k^2\log k = (B/2e^{1+\gamma})L$. Then $\log k \sim \frac12 \log L$,
and we have
$$
m_T \ge e^k\ge  \exp\left(\sqrt {\frac{BL}{e^{1+\gamma}\log L}}
\right).
$$
Note that if $B=2$, then $(B/e^{1+\gamma})^{1/2}=0.64...$ , and
if $B=1$ it equals $0.45...$.
Of course, this assumes
we have uniformity in $k$ out to $\sqrt L$.
On the other hand, if $B=0$ in (29), so that $g_k$
has the form suggested by Keating and Snaith,
then we find that the maximum
is attained when $k\log k$ is near $(A/2e^{1/2+\gamma})L$.
This implies that $\log k \sim \log L$ and that
$$
m_T \ge  \exp\left(\frac{AL}{4e^{\gamma}\log L}\right).
$$
Thus, if $g_k$ grows as suggested by (23), and if (29)
holds uniformly for $k\ll L$, then (24) will be closer
to the truth than (25). Similarly, (26) and (27) would
reflect the true order of $S(T)$ and $|\zeta(1+it)|$,
respectively. Analogously, this suggests that (28)
reflects the correct maximal size of $|L(1,\chi)|$.
This is at odds with the usual view in these questions
which is that the true order is most likely to be near
the lower bounds.
See the forthcoming paper of Granville and Soundararajan [GS]
where one may find similar computations and
and an asymptotic evaluation of
$\sum_{n} d_k(n)^2/n^2$
which resembles the result of the above Proposition.

The second author wishes to express his sincere gratitude
to the American Institute of Mathematics for its generous support
and hospitality while he was working on this paper.

\medskip

\medskip

\centerline{\bf The conjectural formula for $D_k(x,h)$}

\medskip

In this section we sketch the derivation of the form and properties
of $m_k(x,h)$, the conjectural main term for
$$
D_k(x,h)=\sum_{n\le x}d_k(n)d_k(n+h).
$$

We assume that $(a,q)=1$ and that the main term for the sum
$\sum_{n\le x}d_k(n)e(\frac{an}{q})$ can be written in the form
$\frac{1}{q}\int_0^x P_{k}(y,q)~dy$ independently of $a$.
Then applying the $\delta$--method of
Duke, Friedlander, and Iwaniec [DFI] in exactly the same way they do
in the case $k=2$, but ignoring all error
terms, immediately leads to
$$
m'_{k}(x,h)=\sum_{q=1}^{\infty}\frac{c_{q}(h)}{q^{2}}P_{k}(x,q)^{2}\ ,
\tag 30
$$
where $c_{q}(h)=\sum_{d\mid q, d\mid h}d\mu(\frac{q}{d})$ is Ramanujan's sum.
Substituting this expression for $c_{q}(h)$ into the right--hand side,
changing
the order of summation, which will be justified by absolute
convergence once we have established the bound for $P_{k}(x,q)$ in (38 )
below, and relabeling variables, we find that
$$
m_{k}'(x,h)=\sum_{d\mid h} \frac{f_k(x,d)}{d}\ ,
$$
where
$$
f_k(x,d)=\sum_{q=1}^\infty \frac{\mu(q)}{q^2}P_k(x,qd)^2\ .
$$

Next we need to determine an explicit expression for $P_{k}(x,q)$. To
do this we consider the generating function
$$
D_{k}(s,\frac{a}{q})=\sum_{n=1}^{\infty}d_k(n)e(\frac{an}{q})n^{-s}\ ,
\tag 31
$$
with $(a,q)=1$ and $\sigma >1$. Now for any integer $m$ we have
$$
e(\frac{m}{q})=\sum_{d\mid m, d\mid q}\phi(\frac{q}{d})^{-1}
\sum_{\chi (mod \frac{q}{d})}
\tau(\overline{\chi}) \chi (\frac{m}{d})\ ,
$$
where the inner sum is over all characters $\chi$
to the modulus $q/d$ and $\tau(\chi)=\sum_{b (mod q)}
\chi(b)e(\frac{b}{q})$ is Gauss' sum. Using this to replace the
exponential in (31) and rearranging the resulting sums, we obtain
$$
D_{k}(s,\frac{a}{q})=q^{-s}\sum_{d\mid q}\phi(d)^{-1}d^{s}
\sum_{\chi (mod q)}\chi(a)\tau(\overline{\chi})
\sum_{m=1}^{\infty}d_{k}(\frac{qm}{d})\chi(m)m^{-s}\ .
$$
Now if $r=\prod_pp^{\alpha}$, we see that
$$
\eqalign{
\sum_{m=1}^{\infty}d_{k}(rm)\chi(m)m^{-s}=&
\prod_{p\mid
r}\left(\frac{\sum_{j=0}^{\infty}d_{k}(p^{j+\alpha})\chi(p^{j})p^{-js}}
{\sum_{j=0}^{\infty}d_{k}(p^{j})\chi(p^{j})p^{-js}}\right) \cr
&\times
\prod_{p}\left(\sum_{j=0}^{\infty}d_{k}(p^{j})\chi(p^{j})p^{-js}\right) \cr
=&g_{k}(s,r,\chi) L^{k}(s,\chi),
}\tag 32
$$
say. Thus, for $\sigma >1$ we have
$$
D_{k}(s,\frac{a}{q})=q^{-s}\sum_{d\mid q}\phi(d)^{-1}d^{s}
\sum_{\chi (mod q)}\chi(a)\tau(\overline{\chi}) g_{k}(s,q/d,\chi)
L^{k}(s,\chi)\ .
$$
This provides a meromorphic continuation of $D_{k}(s,\frac{a}{q})$ to
the whole complex plane and shows that its only possible pole in $\sigma
>0$ occurs at $s=1$ and is due to the principal character
$\chi_{d}^{(0)} \pmod{d}$ for each $d$ dividing $q$.
Thus the singular part of
$D_{k}(s,\frac{a}{q})$ is the same as that of
$$
q^{-s}\sum_{d\mid q}\phi(d)^{-1}d^{s}
\chi_{d}^{(0)}(a)\tau(\chi_{d}^{(0)})
\sum_{m=1}^{\infty}d_{k}(\frac{qm}{d})\chi_{d}^{(0)}(m)m^{-s}\ .
$$
Now $\tau(\chi_{d}^{(0)})=\sum_{\Sb b (mod\ d)\\ (b,d)=1 \endSb}
e(\frac{b}{d})=c_{d}(1)=\mu(d)$,
and $\chi_{d}^{(0)}(m)=\sum_{e\mid m, e\mid d}\mu(e)$, so we find
that this equals
$$
q^{-s}\sum_{d\mid q}\frac{\mu(d)}{\phi(d)}d^{s}
\sum_{e\mid d}\mu(e)e^{-s}
\sum_{n=1}^{\infty}d_{k}(\frac{qen}{d})n^{-s}\ . \tag 33
$$
We can use (32) to express the sum over $n$ here in terms of the
zeta--function.
Taking $\chi$ equal to $\chi_{1}^{(0)}$, the principal character
$\pmod{1}$, and writing $g_k(s,r)=g_k(s,r,\chi_{1}^{(0)})$,
we deduce from (32) that
$$
\sum_{n=1}^{\infty}d_{k}(rn)n^{-s}=g_{k}(s,r)\zeta^{k}(s),
$$
where
$$
g_{k}(s,r)=\prod_{p\mid r}\left( (1-p^{-s})^{k}\sum_{j=0}^{\infty}
d_{k}(p^{j+\alpha})p^{-js}\right)\ .
$$
Inserting this into (33), we find that the singular part of
$D_{k}(s,\frac{a}{q})$ is identical to that of
$$
q^{-s}\zeta^{k}(s)G_{k}(s,q),
$$
where
$$
G_{k}(s,q)=\sum_{d\mid q}\frac{\mu(d)}{\phi(d)}d^{s}
\sum_{e\mid d}\mu(e)e^{-s}g_{k}(s,qe/d) \ ,
$$
and we note that this is independent of $a$.
From this and Perron's formula we see that $\frac{1}{q}
\int_0^x P_{k}(y,q)~dy$, the main term for
$\sum_{n\le x}d_k(n)e(\frac{an}{q})$,
should be given by
$$
\frac{1}{2\pi i}\int_{|s-1|=1/8}
\zeta(s)^kG_k(s,q)\frac{\left(x/q\right)^{s}}{s}~ds\ .
$$
Thus, differentiating with respect to $x$, we find that
$$
P_k(x,q)=\frac{1}{2\pi i}\int_{|s-1|=1/8}
\zeta(s)^kG_k(s,q)\left(\frac{x}{q}\right)^{s-1}~ds \ .
\tag 34
$$
On changing $s$ to $s+1$, we obtain (15).

It only remains to prove (18). However, before doing this we
derive a formula that we feel is interesting in its
own right, even though we do not require it here. Define
$$
\Cal{D}_{k}(s, h)=\sum_{n=1}^{\infty}d_{k}(n)d_{k}(n+h)n^{-s}.
$$
Then
$$
m_{k}(x, h)=\frac{1}{2\pi i}\int_{|s-1|=1/8}
\Cal{D}_k(s, h)\frac{x^{s}}{s}~ds \ ,
$$
so we see that
$$
m'_{k}(x, h)=\res_{s=1} \Cal{D}_k(s, h)x^{s-1}.
$$
On the other hand, since $D_{k}(s,\frac{1}{q})$ and
$q^{-s}\zeta^{k}(s)G_{k}(s,q)$ have the same singular part at
$s=1$, from (34) we have that
$$
\frac {1}{q}P_{k}(x, q)= \res_{s=1}D_{k}(s,\frac{1}{q})x^{s-1}.
$$
It therefore follows from (30) that
$$
\res_{s=1} \Cal{D}_k(s, h)x^{s-1}=\sum_{q=1}^{\infty}c_{q}(h)
\left(\res_{s=1}D_{k}(s,\frac{1}{q})x^{s-1}\right)^{2}.
$$
This is the formula referred to above.

We now proceed to the proof of (18).
From (16) and (17) we see that $g_{k}(s,1)=G_{k}(s,1)=1$, and that for
$\alpha \ge 1$
$$
G_{k}(s,p^{\alpha})=(1-\frac{1}{p})^{-1}
\left( g_{k}(s,p^{\alpha})-p^{s-1}g_{k}(s,p^{\alpha-1}) \right).
$$
Using this with (17) and the easily proven identity
$$
d_k(p^{\beta})=d_{k-1}(p^{\beta}) + d_k(p^{\beta-1}) \ , \tag 35
$$
we obtain
$$
G_{k}(s, p^{\alpha})=(1-\frac{1}{p})^{-1}(1-\frac{1}{p^{s}})^{k}
\sum_{j=0}^{\infty}p^{-js}\left(d_{k-1}(p^{\alpha+j}) +
(1-p^{s-1})d_k(p^{\alpha+j-1})\right) \ .\tag 36
$$
It is also not difficult to show, for example by induction on $\beta$, that
$$
d_l(p^{\alpha+\beta}) \le d_{l}(p^{\alpha})d_l(p^{\beta}) \ . \tag 37
$$
Therefore we have that
$$
|G_{k}(s,p^{\alpha})| \le (1-\frac{1}{p})^{-1}|1-\frac{1}{p^{s}}|^{k}
\left( (d_{k-1}(p^{\alpha})(1-p^{-\sigma})^{-k+1} +
d_k(p^{\alpha-1})(1-p^{-\sigma})^{-k} |1-p^{s-1}|\right) \ .
$$
Now we assume that $p^{\alpha} \mid \mid q$ and we restrict s to the
circle $|s-1|=\frac{A_{1}}{\log qx}$. Then there exist positive
constants $A_{2}$ and $A_{3}$ such that
$$
 (1-\frac{1}{p}) e^{-\frac {A_{2}\log p}{p\log qx}} \ \le \ \
 |1-\frac{1}{p^{s}}| \ \le \ \
 (1-\frac{1}{p}) e^{\frac {A_{2}\log p}{p\log qx}}
$$
and
$$
 |1-p^{s-1}|\ \le \ A_3 \log p/ \log qx \ .
$$
Hence, for some positive constant $A_{4}$ we have
$$
|G_{k}(s, p^{\alpha})| \le e^{\frac {A_{4}k \log p}{p \log qx}}
\left(d_{k-1}(p^{\alpha}) +
d_{k}(p^{\alpha-1})(1-\frac{1}{p})^{-1} A_3 \log p/ \log qx \right).
$$
We now use the simple formula
$$
d_{k}(p^{\alpha-1})= \frac{\alpha}{k-1}d_{k-1}(p^{\alpha})
$$
and find that
$$
\eqalign{
|G_{k}(s, p^{\alpha})|
&\le d_{k-1}(p^{\alpha})e^{\frac {A_{4}k\log p}{p\log qx}}
\left( 1+A_{5}\alpha \log p/k\log qx \right)  \cr
&\le d_{k-1}(p^{\alpha})e^{\frac {B_{k}\log p^{\alpha}}{\log qx}}\ ,
}
$$
where $B_{k}$ depends only on k.
It follows that
$$
|G_{k}(s, q)| \le  d_{k-1}(q)
e^{\frac {B_{k}\log q}{\log qx}} \ll_{k} d_{k-1}(q).
$$
Since $\zeta(s)\ll |s-1|^{-1}$ near $s=1$, if we use this in (34) and
shrink the
path of integration to the circle $|s-1|=\frac{A_{1}}{\log qx}$, we
deduce that
$$
P_k(x,q)\ll d_{k}(q)\log^{k-1}qx\ . \tag 38
$$
Finally, we use this bound in (14), separate the resulting divisor functions
by means of (37), and recall that $d\le x$ to obtain
$$
f_{k}(x, d)\ll d^{2}_{k-1}(d) \log^{2k-2}x\ ,
$$
which is (18).

\medskip

\medskip

\centerline{\bf The generating function $F_k(x,z)$}

\medskip

In applying Conjecture 3 in the next section, it will turn out that what
we actually require is the behavior of the generating function
$$
F_k(x,z)=\sum_{d=1}^\infty \frac{f_{k}(dx,d)}{d^{z+1}}
$$
near $z=0$. By (14) and (15) we see that
$$
F_k(x,z)=\frac{1}{(2\pi i)^2}\iint
\Sb |s|=1/8 \\|w|=1/8 \endSb
\zeta^{k}(s+1)\zeta^{k}(w+1)\Cal{H}_{k}(z, s, w) x^{s+w}~ds~dw\ ,
\tag 39
$$
where
$$
\Cal{H}_{k}(z, s, w)=\sum_{d=1}^{\infty}\frac{1}{d^{1+z}} \sum_{q=1}^{\infty}
\frac{\mu(q)G_{k}(s+1, dq)G_{k}(w+1, dq)}{q^{2+s+w}}\ .
$$
If we define
$$
\eqalign{
h_{k}(p^{\alpha})=&h_{k}(s, w, p^{\alpha})\cr
=&G_{k}(s+1, p^{\alpha})G_{k}(w+1, p^{\alpha})-
G_{k}(s+1, p^{\alpha+1})G_{k}(w+1, p^{\alpha+1})p^{-2-s-w}\ ,
}\tag 40
$$
then from the multiplicativity of $G_{k}(s+1, q)$ we see that the sum
over $q$ in the definition of $\Cal{H}_{k}(z, s, w)$ equals
$$
\prod_{p\mid \mid d}\left(\frac{h_{k}(p^{\alpha})}{h_{k}(p)}\right)
\prod_{p}h_{k}(p)\ .
$$
The first product defines a multiplicative function of d, so we find
that
$$
\eqalign{
\Cal{H}_{k}(z, s, w)=&\left(\prod_{p}h_{k}(p)\right) \sum_{d=1}^{\infty}
 \frac{1}{d^{1+z}}
\prod_{p\mid \mid d}\left(\frac{h_{k}(p^{\alpha})}{h_{k}(p)}\right) \cr
=&\left(\prod_{p}h_{k}(p)\right) \prod_{p}\left(\sum_{\alpha=0}^{\infty}
\frac {h_{k}(p^{\alpha})/h_{k}(p)}{p^{\alpha(1+z)}}\right) \cr
=&\prod_{p}\left(\sum_{\alpha=0}^{\infty}
\frac { h_{k}(p^{\alpha})}{p^{\alpha(1+z)}}\right).
} \tag 41
$$
Now $G_{k}(s+1,1)=1$, and for $\alpha \ge 1$
$$
G_{k}(s+1, p^{\alpha})=d_{k-1}(p^{\alpha})+d_{k}(p^{\alpha-1})(1-p^{s})
+O_{k}(p^{-1-\sigma})+O_{k}(p^{-1})
$$
by (36). Thus, writing $\Re s=\sigma$, $\Re w=u$, and $\Re z=x$, we see that
$$
h_{k}(1)=1+O_{k}(p^{-2+|\sigma|+|u|})\ ,
$$
$$
\eqalign{
h_{k}(p^{\alpha})=&\left(d_{k-1}(p^{\alpha})+d_{k}(p^{\alpha-1})(1-p^{s})
\right)
\left(d_{k-1}(p^{\alpha})+d_{k}(p^{\alpha-1})(1-p^{w})\right) \cr
&+O_{k}(p^{-1+|\sigma|+|u|})\ ,
}
$$
and
$$
\sum_{\alpha=0}^{\infty}
\frac { h_{k}(p^{\alpha})}{p^{\alpha(1+z)}}=
1+\frac {(k-p^{s})(k-p^{w})}{p^{1+z}}
+O_{k}(p^{-2+|\sigma|+|u|+|x|}) \ ,
$$
where the function bounded by the $O$--term in the last line
is analytic in the region $|\sigma|+|u|+|x| <2$. When we use this with
(41)
we obtain
$$
\eqalign{
\Cal{H}_{k}(z, s, w)
=&\prod_{p}\left(1+\frac{k^{2}}{p^{1+z}}-\frac{k}{p^{1+z-s}}-\frac{k}{p^{1+z-w}}
+\frac{1}{p^{1+z-s-w}}+ \dots \right) \cr
=&\zeta^{k^{2}}(1+z)\zeta^{-k}(1+z-s)\zeta^{-k}(1+z-w)
\zeta(1+z-s-w)\Cal{H}^{*}_{k}(z, s, w) \ ,
}
$$
where
$$
\eqalign{
\Cal{H}^{*}_{k}&(z, s, w)\cr
&=\prod_{p}\left(
 (1-\frac{1}{p^{1+z}})^{k^{2}}
(1-\frac{1}{p^{1+z-s}})^{-k}
(1-\frac{1}{p^{1+z-w}})^{-k}
(1-\frac{1}{p^{1+z-s-w}})
(\sum_{\alpha=0}^{\infty}\frac { h_{k}(p^{\alpha})}{p^{\alpha(1+z)}})
\right)
}
$$
is analytic for $|\sigma|+|u|+|x| <1$. Combining this with (39), we find that
$$
\eqalign{
F_k&(x,z) \cr
&=\frac{1}{(2\pi i)^2}\iint \Sb |s|=1/8 \\|w|=1/8 \endSb
\frac{\zeta(s+1)^k\zeta(w+1)^k\zeta(z+1)^{k^2}\zeta(z+1-s-w)x^{s+w}}
{\zeta(z+1-s)^k\zeta(z+1-w)^k}\Cal{H}^{*}_{k}(z, s, w)~ds~dw \ .
} \tag 42
$$

Finally we evaluate $\Cal{H}^{*}_{k}(0, 0, 0)$ as this is also required
in the next section. By (40) and the definition of $\Cal{H}^{*}_{k}(z, s,
w)$ we have
$$
\Cal{H}^{*}_{k}(0, 0, 0)=\prod_{p}\left((1-\frac{1}{p})^{(k-1)^{2}}
\sum_{\alpha=0}^{\infty}\frac{(G^{2}_{k}(1, p^{\alpha})
-G^{2}_{k}(1, p^{\alpha+1})p^{-2})}{p^{\alpha}}
\right)\ . \tag 43
$$
By (36), $G_{k}(1, p^{\alpha})=(1-\frac{1}{p})^{k-1}
\sum_{j=0}^{\infty}d_{k-1}(p^{\alpha+j})p^{-j}$.
Hence, factoring the numerator in the sum over $\alpha$ as a diference
of two squares, we obtain
$$
\eqalign{
G^{2}_{k}(1, p^{\alpha})&
-G^{2}_{k}(1, p^{\alpha+1})p^{-2} \cr
&=(1-\frac{1}{p})^{2k-2}d_{k-1}(p^{\alpha})
\left(d_{k-1}(p^{\alpha})+2\sum_{j=0}^{\infty}
\frac{d_{k-1}(p^{\alpha+j+1})}{p^{j+1}}\right ).
}
$$
We therefore see that the typical factor in (43) equals
$$
(1-\frac{1}{p})^{k^{2}-1}
\left(
\sum_{\alpha=0}^{\infty}
\frac{d^{2}_{k-1}(p^{\alpha})}{p^{\alpha}}+
2\sum_{\alpha=1}^{\infty}
\frac{d_{k-1}(p^{\alpha})}{p^{\alpha}}\left
(\sum_{l=0}^{\alpha-1}d_{k-1}(p^{l})\right )
\right)\ .
$$
Since $\sum_{q\mid n}d_{k-1}(q)=d_{k}(n)$, the sum over $l$ equals
$d_{k}(p^{\alpha-1})$, so this is
$$
(1-\frac{1}{p})^{k^{2}-1}
\left(1+
\sum_{\alpha=1}^{\infty}
\frac{\left (d_{k-1}(p^{\alpha})+d_{k}(p^{\alpha-1}\right
)^{2}-d^{2}_{k}(p^{\alpha-1})}{p^{\alpha}}
\right).
$$
We now use (35) to see that this is
$$
\eqalign{
&(1-\frac{1}{p})^{k^{2}-1}
\left(
\sum_{\alpha=0}^{\infty}
\frac{d^{2}_{k}(p^{\alpha})}{p^{\alpha}}-
\sum_{\alpha=1}^{\infty}
\frac{d^{2}_{k}(p^{\alpha-1})}{p^{\alpha}}
\right)\cr
\ \ \ &=(1-\frac{1}{p})^{k^{2}}
\left(
\sum_{\alpha=0}^{\infty}
\frac{d^{2}_{k}(p^{\alpha})}{p^{\alpha}}
\right).
}
$$
Inserting this into (43) and comparing with (10), we deduce that
$$
\Cal{H}^{*}_{k}(0, 0, 0)=a_{k}\ .
$$

\medskip

\medskip

\centerline{\bf Conjecture 4}

\medskip

A precise version of the mean value formula we require to derive the
formula in Conjecture 4 is given in the paper of Goldston and Gonek [GG].
However, we adopt a simpler heuristic approach that will lead to the
same formula more quickly.

Recall that
$$
\D_{k,N}(s)=\sum_{n=1}^N\frac {d_k(n)}{n^s}
$$
and that we are to estimate
$$
\Cal{I}(T)=\Cal{I}_{k,N}(T)=\int_T^{2T}|\D_{k,N}(1/2+it)|^2~dt.
$$
The reason we integrate from $T$ to $2T$ rather than from $0$ to $T$
is that for $t$ near $0$ the integrand is very large, on the order of
$N^{1/2}$, so the mean square for small $t$ is about $N$. Since $N$
can be as large as $T^2$, this would dominate, and so obscure, the behavior
of the mean value away from the real axis. We let
$$
\Cal{J}(T)=\Cal{J}_{k,N}(T)=\frac 12 \int_{-T}^T |\D_{k,N}(1/2+it)|^2~dt
\tag 44
$$
and then obtain our estimate for $\Cal{I}(T)$ via the formula
$$
\Cal{I}(T)=\Cal{J}(2T)-\Cal{J}(T)\ . \tag 45
$$
Squaring out and integrating term--by--term in (44), we find that
$$
\split
\Cal{J}(T)&=T\sum_{n=1}^N \frac {d_k(n)^2}{n} + \sum_{m\ne
n}\frac{d_k(m)d_k(n)}{(mn)^{1/2}}
\frac {\sin\left(T\log \frac m n\right)}{\log \frac m n }\\
&=\Cal{J}_d(T)+\Cal{J}_o(T)\ ,
\endsplit \tag 46
$$
say. Using the symmetry in $m$ and $n$ and writing $m=n+h$, we find that
$$
\Cal{J}_o=2\sum_{n=1}^N\sum_{h=1}^{N-n}\frac{d_k(n)d_k(n+h)}{n\sqrt{1+h/n}}
\frac{\sin\left(T\log\left(1+\frac h n\right)\right)}
{\log\left(1+\frac {h}{n} \right)}\ .
$$

Now we make several approximations which are justified in the paper of
Goldston and Gonek for the weighted version of this formula.
Namely, we replace $\log(1+h/n)$
by $h/n$, $\sqrt{1+h/n}$ by 1, and
$d_k(n)d_k(n+h)$ by $m_k'(n,h)$ from Conjecture 3.
Finally, the sum over $n$ can be replaced by an integral and the
sum over $h$ extended to $\infty$.
This leads to
$$
\Cal{J}_o\sim 2\int_0^N \sum_{h=1}^\infty
\frac{m_{k}'(x,h)}{h}\sin\left(\frac{Th}{x}\right)~dx\ .
$$
By (13), the right--hand side equals
$$
\split
&2\int_0^N\sum_{h=1}^\infty \frac{1}{h}\sum_{d\mid h}
\frac{f_{k}(x,d)}{d}\sin\left(\frac{Th}{x}\right)~dx
\\
=&2\int_0^N\sum_{d=1}^\infty\frac{f_{k}(x,d)}{d^2}\sum_{h=1}^\infty
\frac{\sin(Thd/x)}{h}~dx
\\
=&-2\pi\int_0^N\sum_{d=1}^\infty\frac{f_{k}(x,d)}{d^2}\left(
\left\{\frac{Td}{2\pi x}\right\}-\frac12\right)~dx \ ,\endsplit
$$
where
$$\{x\}=\cases x-[x] \text{\ if $x$ is not an integer,}\\
1/2 \text{\ if $x$ is an integer}\ \ . \endcases
$$
The ``-1/2'' term  leads to a large contribution if
$N$ is large, but it is independent of $T$ and so disappears when
we take the difference $\Cal{J}(2T)-\Cal{J}(T)=\Cal{I}(T)$.
Thus, we may express $\Cal{J}_o$ as
$$
\Cal{J}_o(T)=\Cal{J}_{o,1}(T)+\Cal{J}_{o,2}\ , \tag 47
$$
where $\Cal{J}_{o,1}(T)$ is the part with
``$\{\}$'' and $\Cal{J}_{o,2}$ is the part with ``-1/2''.

We now make the change of variable $y=Td/(2\pi x)$ in $\Cal{J}_{o,1}$
and find that
$$
\eqalign {
\Cal{J}_{o,1}(T)&=-2\pi\int_0^N\sum_{d=1}^\infty\frac{f_{k}(x,d)}{d^2}
\left\{\frac{Td}{2\pi x}\right\}~dx \cr
&=-T\sum_{d=1}^\infty \frac 1 d \int_{\frac{Td}{2\pi N}}^\infty
\frac{f_{k}(Td/2\pi y,d)}{y^2}\{y\}~dy.
}
$$
We split the sum over $d$ into $d\le 2\pi N/T$ and $d> 2\pi N/T$.
From our estimate for $f_{k}$ in (18) we find that the contribution
to $\Cal{J}_{o,1}$ from the upper range of $d$ is
$$
\ll T\sum_{N/T\ll d}\frac{\tau_{k-1}(d)^2}{d}\frac{N}{Td}L^{2k-2}
\ll TL^{(k-1)^2-1+2k-2}=TL^{k^2-2}.
$$
In the lower range of $d$ we split the integral over $y$ into two ranges:
$Td/2\pi N\le y<1$ and $y\ge 1$.
The contribution from the upper range of $y$ is
$$
\ll T\sum_{d\ll N/T}\frac{\tau_{k-1}(d)^2}{d}L^{2k-2}\frac{N}{dT}\ll
TL^{k^2-1}.
$$
Hence, since $\{y\}=y$ for $0<y<1$, we see that
$$
\eqalign{
\Cal{J}_{o,1}(T)
&=-T\sum_{d\le N/T} \frac 1 d \int_{\frac{Td}{2\pi N}}^1
\frac{f_{k}(Td/2\pi y,d)}{y}~dy +O(TL^{k^2-1}) \cr
&=-T\int_{\frac{T}{2\pi N}}^1\sum_{d\le yN/T}
\frac{f_{k}(Td/2\pi y,d)}{d}\frac{dy}{y} +O(TL^{k^2-1}) \ .
}
$$

We now use the expression in (42) for the generating function
$F_{k}(x,z)$ of $f_{k}(dx,d)$.
In doing this, we retain only the first term in the Laurent
expansion at zero of the various factors in the integrand
and find that
$$
\Cal{J}_{o,1}(T)
\sim -a_kT\int_\frac{T}{N}^1\frac{1}{(2\pi i)^3}\iiint
\Sb |s|=\frac {1}{8}\\
|w|=\frac{1}{8}\\ |z|=\frac 12 \endSb
\frac{\left(\frac Ty\right)^{s+w-z}N^z(z-s)^k(z-w)^k}
{z^{k^2+1}(z-s-w)s^kw^k}~ds~dw~dz\frac{dy}{y}
$$
To evaluate this we make the substitutions $s\to sz$
and $w\to wz$, and carry out the integration in $z$. We then find
that
$$
\eqalign{
&\Cal{J}_{o,1}(T) \sim \cr
&-\frac{a_k T}{\Gamma(k^2)(2\pi i)^2 }\int_\frac{T}{N}
\iint \Sb |s|=\frac {1}{16}\\
|w|=\frac{1}{16}  \endSb
\frac{\left((s+w-1)\log \frac Ty +\log N\right)^{k^2-1}
(1-s)^k(1-w)^k}{(1-s-w)s^kw^k}~ds~dw\frac{dy}{y}.
} \tag 48
$$

Next we write
$$
N=T^{1+\eta} \ ,
$$
where $0\le \eta \le 1$, and make the substitution
$$
y=T^{\alpha(1+\eta)-\eta}.
$$
Then $y=1$ corresponds to $\alpha=1-1/(1+\eta)$ and $y=T/N$ corresponds
to $\alpha=0$. Also,
$$
(s+w-1)\log\frac Ty+\log N=\log N(1+(s+w-1)(1-\alpha))
$$
 and
$$
\frac{dy}{y}=\log N d\alpha.
$$
Thus the asymptotic expression for (48) becomes
$$
\Cal{J}_{o,1}(T)\sim
\frac{a_k}{\Gamma(k^2)}(1+\eta)^{k^2} TL^{k^2}\Cal{M}_k\ ,\tag 49
$$
where
$$
\Cal{M}_k=-\frac{1}{(2\pi i)^2}
\int_0^{1-\frac 1{1+\eta}}\iint
 \Sb |s|=\frac{1}{16}\\
|w|=\frac{1}{16} \endSb
\frac{(1+(s+w-1)(1-\alpha))^{k^2-1}
(1-s)^k(1-w)^k}{(1-s-w)s^kw^k}~ds~dw~d\alpha \ .
$$

We observe for future reference that $\Cal{M}_k\to 0$ as $k \to \infty$.
In fact, since
$$
1+(s+w-1)(1-\alpha)=\alpha+(s+w)(1-\alpha),
$$
we find that
$$
|\Cal{M}_k|\le (1-1/(\eta+1)+1/8)^{k^2-1}{15}^{2k}\ll (4/5)^{k^2}
$$
as $k\to \infty$.

Next, we expand $(1+(s+w-1)(1-\alpha))^{k^2-1}$ into powers of
$(s+w-1)$ and find that
$$
(1-(1-s-w)(1-\alpha))^{k^2-1}=\sum_{n=0}^{k^2-1}\binom{k^2-1}{n}
(-1)^n(1-s-w)^n(1-\alpha)^n \ .
$$
Thus
$$
\Cal{M}_k=-\int_0^{1-\frac 1{1+\eta}}\sum_{n=0}^{k^2-1}(-1)^n
\gamma_k(n)(1-\alpha)^n~d\alpha \ ,
$$
where
$$
\gamma_k(n)=
\frac{1}{(2\pi i)^2}\iint \Sb |s|=\frac{1}{16}\\
|w|=\frac{1}{16}\endSb
\frac{(1-s-w)^{n-1} (1-s)^k(1-w)^k}{s^kw^k}~ds~dw\ .
$$
Now
$$\split
\gamma_k(0)&=\frac{1}{(2\pi i)^2}\iint \Sb |s|=\frac{1}{16}\\
|w|=\frac{1}{16}\endSb
\frac{(1-s)^k(1-w)^k}{(1-s-w)s^kw^k}~ds~dw\\
&=\sum_{i=0}^k\sum_{j=0}^k\binom{k}{i}\binom{k}{j}(-1)^{i+j}
\frac{1}{(2\pi i)^2}\iint \Sb |s|=\frac{1}{16}\\
|w|=\frac{1}{16} \endSb
\frac{~ds ~dw }{(1-s-w)s^iw^j}~ds~dw\\
&=\sum_{i=0}^k\sum_{j=0}^k\sum_{m=0}^\infty\binom{k}{i}\binom{k}{j}(-1)^{i+j}
\frac{1}{(2\pi i)^2}\iint \Sb |s|=\frac{1}{16}\\
|w|=\frac{1}{16}\endSb
\frac{(s+w)^m ~ds ~dw }{s^iw^j}~ds~dw\\
&=\sum_{i=0}^k\sum_{j=0}^k\binom{k}{i}\binom{k}{j}\binom{i+j-2}{i-1}(-1)^{i+j}.
\endsplit
$$
Actually, this can be simplified to $\gamma_k(0)=k$, but this fact
is not necessary to proceed.
In a similar manner we find that
$$
\gamma_k(n)=(-1)^n\sum_{i=0}^k\sum_{j=0}^k\binom{k}{i}\binom{k}{j}
\binom{n-1}{i-1,j-1,n-i-j+1}\ .
$$
Finally, integrating with respect to $\alpha$, we
obtain
$$
\Cal{M}_k=-\sum_{n=0}^{k^2-1}(-1)^n\gamma_k(n)\frac{\left(1-(\eta+1)^{-(n+1)}
\right)}{n+1}.
$$

Combining this, (47), and (49), we obtain
$$
\eqalign{
\Cal{J}_{o}(2T) -\Cal{J}_{o}(T)& =
\Cal{J}_{o,1}(2T) -\Cal{J}_{o,1}(T) \cr
&\sim
-\frac{a_k}{\Gamma(k^2)}(1+\eta)^{k^2} TL^{k^2}
\sum_{n=0}^{k^2-1}(-1)^n\gamma_k(n)\frac{\left(1-(\eta+1)^{-(n+1)}
\right)}{n+1}\ .
} \tag 50
$$
Also, from (9) we have
$$
\Cal{J}_d(2T)-\Cal{J}_d(T)
\sim \frac{a_k}{\Gamma(k^2+1)}(1+\eta)^{k^2}L^{k^2} \ .
$$
It now follows from (45), (46), (50), and this that
$$
\eqalign{
\Cal{I}_{k,N}&=\int_T^{2T}\left|\D_{k,N}(1/2+it)\right|^2~dt \cr
&\sim w_{k}(\eta)\frac{a_k}{\Gamma(k^2+1)}TL^{k^2}\ ,
}
$$
where
$$
w_{k}(\eta)=(1+\eta)^{k^2}\left(1-\sum_{n=0}^{k^2-1}\binom{k^2}{n+1}
\gamma_k(n)\frac{\left(1-(1+\eta)^{-(n+1)}\right)}{n+1}\right) \ .
$$
This is Conjecture 4.

\medskip

\medskip

\centerline{\bf The sixth and eighth power moment conjectures}

\medskip

We first note an alternative expression for $a_{k}$.
Conrey and Ghosh [CG3] have shown that $a_{k}=a_{1-k}$. Thus, by (10),
$$
a_k=\prod_p\left(\left(1-\frac 1 p\right)^{(k-1)^2}\sum_{r=0}^{\infty}
\frac{d^{2}_{1-k}(p^{r})}{p^r}\right).
$$
Since $d_k(p^r)=\binom{k+r-1}{r}=(-1)^{r}\binom{-k}{r}$, we see that
$d_{1-k}(p^{r})=(-1)^{r}\binom{k-1}{r}$, so we have that
$$
a_k=\prod_p\left(\left(1-\frac 1 p\right)^{(k-1)^2}\sum_{r=0}^{k-1}
\frac{{\binom{k-1}{r}}^{2}}{p^r}\right)\ . \tag 51
$$
This is the expression for $a_{k}$ we have used in (3) and Conjecture
1.

By (6) we expect that
$$
\int_T^{2T} |\zeta(1/2+it)|^6~dt \sim
\int_T^{2T}\left(\left|\D_{3,T^{1+ \eta}}(1/2+it)\right|^2+
\left|\D_{3,T^{2-\eta}}(1/2+it)\right|^2\right)~dt
$$
for any $\eta$ with $0\le \eta \le 2$. Using Conjecture 4 and adding the
results together for $T/2, T/4, \dots $, we obtain
$$
\eqalign{
I_3(T)=\int_0^{T} |\zeta(1/2+it)|^6~dt
\sim (w_{3}(\eta)+w_{3}(1-\eta))\frac{a_3}{9!}TL^9 \ .
}
$$
From the formula for $w_{k}(\eta)$ in the conjecture we calculate that
$$
w_{3}(\eta)=1+9\eta+36\eta^2+84\eta^3+126\eta^4
-630\eta^5+588\eta^6+180\eta^7-9\eta^8+2\eta^9 \ ,
$$
and it is not difficult to verify that
$$
w_{3}(\eta)+w_{3}(1-\eta)=42
$$
for all $\eta$. This identity provides compelling evidence for the
sixth power moment conjecture of Conrey and Ghosh.

Similarly we calculate that
$$
\split
w_{4}(\eta)=&1+16\eta+120\eta^2+560\eta^3+1820\eta^4+4368\eta^5+8008\eta^6+11440
\eta^7\\
&+12870\eta^8+
11440\eta^9
-152152\eta^{10}+179088\eta^{11}-78260\eta^{12}\\
&+14000\eta^{13}-1320\eta^{14}+16\eta^{15}-3\eta^{16}\ .
\endsplit
$$
We then find that
$$
w_{4}(1)=12012 \ ,
$$
and this leads to
$$
I_4(T)\sim 24024\frac{a_4}{16!}TL^{16} \ ,
$$
which is Conjecture 1.

\medskip

\medskip

\centerline{\bf Proof of the Proposition}

\medskip

To prove the proposition we work from the expression
for $a_{k}$ given in (51). Using this,
we first prove an upper bound for $a_k$.
Write
$$
a_k=\Pi^-\Pi^{+} \ \ ,
$$
where $\Pi^-$ is the part of the product over the primes
$\le 2k^2$ and $\Pi^+$ is the part over the primes $> 2k^2$.
Then
$$
\split
\Pi^-&\le \prod_{p\le 2k^2}\left(1-\frac1p\right)^{(k-1)^2}
\left(1+\frac {k-1}{p^{1/2}}+\frac{\binom{k-1}{2}}{p}+
\frac{\binom{k-1}{3}}{p^{3/2}}+\dots\right)^2\\
& \le \left(\frac{e^{-\gamma}}{\log
2k^2}(1+o(1))\right)^{(k-1)^2}\prod_{p\le 2k^2}\left(1+\frac{1}{p^{1/2}}\right)^{2(k-1)}\\
&\le \left(\frac{e^{-\gamma}}{2\log k}
\right)^{k^2}e^{o\left(k^2\right)} exp\left(\sum_{p\le
2k^2}\frac{2(k-1)}{p^{1/2}}\right)\\
&=\left(\frac{e^{-\gamma}}{2\log k}\right)^{k^2}e^{o\left(k^2\right)}
e^{O(\left(\frac{k^2}{\log k}\right)}\\
&=\left(\frac{e^{-\gamma}}{2\log k}\right)^{k^2}e^{o\left(k^2\right)}\ .
\endsplit
$$
Next, it is easy to show that
$$
\binom{k-1}{r}^2\le \binom{(k-1)^2}{r}
$$
for $r=0,1,2,\dots$ and $k=1,2,3, \dots$\ , so we see that
$$
\eqalign{
\Pi^{+} \le&\prod_{p> 2k^2}\left(1-\frac 1 p\right)^{(k-1)^2}
\left(1+\frac 1 p\right)^{(k-1)^2} \cr
=&\prod_{p> 2k^2}\left(1-\frac 1 {p^2}\right)^{(k-1)^2}\cr
\le&1 \ .
}
$$
Thus we find that
$$
a_k\le \left(\frac{e^{-\gamma}}{2\log
k}\right)^{k^2}e^{o\left(k^2\right)}\ .
$$

Now we deduce a lower bound. First we have
$$
\eqalign{
\Pi^{-}&
\ge\prod_{p\le 2k^2}\left(1-\frac 1p\right)^{(k-1)^2}\cr
&=\left(\frac{e^{-\gamma}}{\log 2k^2}(1+o(1))\right)^{(k-1)^2}\cr
&=\left(\frac{e^{-\gamma}}{2\log k}\right)^{k^2}e^{o(k^2)}.
}
$$
Also, since $(1-x)^n\ge 1-nx$ for $0<x<1$, we have
$$
\eqalign{
\Pi^+&\ge \prod_{p>2k^2}\left(1-\frac 1
p\right)^{(k-1)^2}\left(1+\frac{(k-1)^2}{p}\right) \cr
&\ge \prod_{p>2k^2}\left(1-\frac {(k-1)^2}p\right)
\left(1+\frac{(k-1)^2}{p}\right) \cr
&\ge  \prod_{p>2k^2}\left(1-\frac {(k-1)^4}{p^2}\right) \cr
&=\exp\left(\sum_{p>2k^2}\log\left(1-\frac{(k-1)^4}{p^2}\right)\right).
}
$$
Since $\log(1-x) \ge -2 x$ for $0\le x\le .8\ ,$ and
$\frac{(k-1)^4}{p^2}\le 0.8$
for $p>2k^2$,
this is
$$
\ge \exp\left(-2\sum_{p>2k^2}\frac{(k-1)^4}{p^2}\right)
\ge \exp\left(-O\left(\frac{k^4}{k^2\log k}\right)\right)=e^{o(k^2)}.
$$
Thus, we find that
$$
a_k\ge \left(\frac{e^{-\gamma}}{2\log
k}\right)^{k^2}e^{o\left(k^2\right)} \ .
$$
Since the upper and lower bounds are the
same, the Proposition follows.

\enddocument